\documentclass[12pt]{article}%
\usepackage{amssymb}
\usepackage{amsfonts}
\usepackage{amsmath}
\usepackage{graphicx}%
\begin{document}

\begin{center}

\textbf{The meet operation in the imbalance lattice of maximal instantaneous
codes: alternative proof of existence}$\footnotetext{
\par
{}}\footnotetext{AMS Classification 2010: \ Primary: 06A07, 94A45}$

\bigskip

Stephan Foldes, D. Stott Parker, S\'{a}ndor Radeleczki

\bigskip

\bigskip

\bigskip\textbf{Abstract }
\end{center}

\textit{An alternative proof is given of the existence of greatest lower
bounds in the imbalance order of binary maximal instantaneous codes of a given
size. These codes are viewed as maximal antichains of a given size in the
infinite binary tree of 0-1 words. The proof proposed makes use of a single
balancing operation instead of expansion and contraction as in the original
proof of the existence of glb.}

\bigskip\bigskip

\textbf{1} \textbf{\ Terminology of codes and introduction}

\bigskip

The set $\left\{  0,1\right\}  ^{\ast}$ of all finite sequences (words) of the
symbols $0$ and $1$ is partially ordered by the \textit{prefix order}
$\leq_{pref}$\ defined by
\[
v\leq_{pref}w\Leftrightarrow\exists z\text{ }vz=w
\]
The prefix-ordered set of words is an infinite binary tree having the empty
word as root. \textit{Instantaneous code}s are defined as the finite
antichains in this tree. (This finiteness shall be assumed throughout the
paper, thus excluding infinite prefix-free sets.) By \textit{lexicographic
(lex) order} we mean the (only) linear extension of the prefix order in which
words incomparable in the prefix order are compared by the \textit{"telephone
book principle"}, i.e. $v0x$ always precedes (is smaller than) $v1y$. We call
an instantaneous code \textit{lex monotone} if the sequence of lengths of the
codewords taken in lex order is monotone (non-decreasing). With respect to a
given lex monotone instantaneous code $C$ whose codewords in lex order are
$c_{1},c_{2},...$, each codeword $c_{i}$ is then identified by its
\textit{index} $i$.

It is well known that for every maximal instantaneous code there is a unique
lex monotone maximal instantaneous code with the same multiset of codeword
lengths. (See e.g. $\left[  FS\right]  $ for more general statements.) Thus
the multiset of codeword lengths, displayed as a list of numbers with
repetitions in non-decreasing order, can be used to denote the maximal
instantaneous lex monotone code, these are the \textit{path-length sequences}
appearing in $\left[  SPPR\right]  $. For example, the only lex monotone
maximal instantaneous code of size $3$ is $\left\{  0,10,11\right\}  $, and
its path-length sequence is $(1,2,2)$. The code can also be displayed by the
binary tree of all prefixes of the codewords (called a \textit{canonical tree
}by Elsholtz, Heuberger and Prodinger in $\left[  EHP\right]  ,$ under
assumption of lex monotonicity), and the path-length sequence is then the
sequence of lengths of root-to-leaf paths of this tree.

For any word $w$ we can use the simplified notation
\[
2^{-w}%
\]
to denote the number obtained by raising $1/2$ to a power equal to the length
of $w$. The \textit{Kraft sum} of any instantaneous code $C$ is then the sum
\[
K(C)=\sum_{w\in C}2^{-w}%
\]
The Kraft sum is always at most $1$, and it is equal to $1$ if and only if the
instantaneous code $C$ is maximal (Kraft $\left[  K\right]  $).

\bigskip

\bigskip

\textbf{2 \ Statements and proofs\bigskip}

With the above terminology and notation we rephrase the definition of
imbalance order and the result that it is a lattice (Stott Parker and Prasad
Ram $\left[  SPPR\right]  $) as follows.

\bigskip

\textbf{Definition of Imbalance Order by Majorization} $\left[  SPPR\right]  $
\textit{Let} $L$ \textit{be the set of lex monotone maximal instantaneous
codes of a same given size. For codes }$A,B\in L,$ \textit{lexicographically
enumerated }%

\[
A=a_{1},a_{2},...\text{ \ \ \ \ \ \ \ }B=b_{1},b_{2},...\text{
\ \ \ \ \ \ \ \ \ \ \ \ \ \ \ \ \ \ \ \ \ \ \ \ \ \ \ \ (1)}%
\]
$A$ \textit{is said to be more balanced (less imbalanced) than or equal to
}$B$,\textit{ in symbols} $A\trianglelefteq B$, \textit{if for all }$m\geq1$
\textit{we have the following inequality for the partial Kraft sums:}%
\[
K(a_{1},...,a_{m})\leq K(b_{1},...,b_{m})
\]

\bigskip

A characterization of the imbalance order via ternary exchanges was also given
in $\left[  SPPR\right]  ,$ and in the sequel we shall give another
characterization by comparing indices in the enumerations (1).

The size $f(t)$ of the poset of lex monotone maximal instantaneous codes of
size $t$ is an exponentially growing function of the parameter $t$, which has
been the object of combinatorial studies since the 1960's (see $\left[
EHP\right]  $, where further references are also given). While a closed
formula for $f(t)$ is not available, Elsholtz, Heuberger and Prodinger
$\left[  EHP\right]  $ gave a new and very tight asymptotic estimate of $f(t)$.

\bigskip

\textbf{Lattice Property of Imbalance Order }$\left[  SPPR\right]  $
\textit{The imbalance-ordered set of lex monotone maximal instantaneous codes
of the same given size is a lattice.}

\bigskip

Due to the lattice property the construction of optimal codes becomes an
optimization problem on a lattice. In a context very different from that of
binary codes, the balance concept introduced in $\left[  SPPR\right]  $ has
also been shown by O'Keefe, Pajoohesh and Schellekens $\left[  OKPS\right]  $
to be relevant in studying the efficiency of algorithms that involve a
bifurcation at each step, as a root-to-leaf path in the decision tree of an
algorithm corresponds to the succession of steps of the algorithm on a
particular input, and path-length corresponds to running time on that input.
Besides pointing to the analogy between the concepts of average codeword
length and average running time, $\left[  OKPS\right]  $ also shows that, with
the exception of the very small lattices, the imbalance lattices are not
modular. Pajoohesh $\left[  P\right]  $ characterizes the most balanced and
the most imbalanced trees in terms of the semilattice structure of the trees themselves.

\bigskip\

In the present paper an alternative proof of the above lattice property is
given, based not on induction on the common size of the codes in the
imbalance-ordered set of codes, but on applying the abstract Criterion below
for a poset to be a lattice. An earlier alternative explanation of the lattice
property, in fact closer to the techniques of the original proof of the result
in $\left[  SPPR\right]  ,$ was given by two of the present authors in
$\left[  FR\right]  $.

\bigskip

\textbf{Criterion for Lattice Property }\textit{For any finite partially
ordered set with minimum and maximum the following conditions are equivalent:
}

\textit{(i) the poset is a lattice,}

\textit{(ii) for every pair of distinct elements} $b,c$, \textit{one of them -
say} $c$ - \textit{can be replaced by a lesser element} \ $d<c$, \textit{such
that} $b$ \textit{and} $c$ \textit{have the same common lower bounds as} $b$
\textit{and} $d,$

\textit{(iii)} \textit{for every pair of distinct elements} $b,c$
\textit{there is an element} $d$ \textit{that is less than} $b$ \textit{or
}$c$, \textit{and such that} $b,c,d$ \textit{have the same common lower bounds
as} $b,c$.

\bigskip

\textbf{Proof} \ Condition (ii) obviously holds in any lattice, while the
existence of a greatest lower bound of $b$ and $c$ is obtained, using
condition (ii), by induction on the number of elements that are below at least
one of $b$ and $c$. Condition (iii) is a re-phrasing of (ii). \ $\square$

\bigskip

The characterization of the imbalance order given below, by comparing indices,
and a reduction lemma, shall make the above Criterion applicable. The
characterization is based on the following description of the comparabilities
between elements of any two finite maximal instantaneous codes:

\bigskip

\textbf{Interval Decomposition Lemma for Two Codes }\textit{For any two
maximal instantaneous codes }$A$ \textit{and} $B$ \textit{there is a unique
positive integer }$n$ \textit{and unique partitions of the lexicographically
ordered codes into }$n$\textit{\ pairwise disjoint non-empty intervals
consecutive in the lexicographic order}%
\begin{align*}
A  &  =A_{1}\cup...\cup A_{n}\text{ \ \ \ \ \ \ \ \ }B=B_{1}\cup...\cup
B_{n}\text{
\ \ \ \ \ \ \ \ \ \ \ \ \ \ \ \ \ \ \ \ \ \ \ \ \ \ \ \ \ \ \ \ \ \ \ (2)}\\
\text{ \ \ \ }A_{1}  &  <...<A_{n}\text{ \ \ \ \ \ \ \ \ }B_{1}<...<B_{n}%
\end{align*}
\textit{such that any words }$x\in A_{i}$ \textit{and }$y\in B_{j}$
\textit{are comparable} \textit{in the prefix order if and only if }$i=j.$

\bigskip

\textbf{Proof }Two elements of $A$ belong to the same interval if and only if
there is some element of $B$ comparable with both. The intervals of $B$ are
defined similarly. The fact that this defines interval decompositions of the
two codes with the same number $n$ of intervals is verified without
difficulty. The claimed properties and uniqueness are also straightforward.
\ $\square$

\bigskip

The interval decompositions (2) also have the following properties:\medskip

(i) for every $i$, at least one of $A_{i}$ or $B_{i}$ is a singleton, both are
singletons if and only if $A_{i}=B_{i}$, otherwise they are disjoint,

(ii) if the interval $A_{i}$ (respectively $B_{i}$) is a singleton, then its
unique element is a prefix of the words in $B_{i}$ (respectively in $A_{i}$),

(iii) for every $i$ we have the equality of the corresponding interval Kraft
sums, $K(A_{i})=K(B_{i})$,

(iv) if $i<j$, $x\in A_{i}\cup B_{i}$, $y\in A_{j}\cup B_{j}$, then $x$ and
$y$ are incomparable in the prefix order and $x$ precedes $y$ lexicographically.

\bigskip

With a view of referring to these interval decompositions in the sequel, we
call the intervals $A_{i}$ (respectively $B_{i}$) in (2) the
\textit{(comparability) blocks} of $A$ with respect to $B$ (of $B$ with
respect to $A$). A block $A_{i}$ (respectively $B_{i}$) is said to be
\textit{dominating} if it is a singleton but $B_{i}$ (respectively $A_{i}$) is
not. In that case the sole element of $A_{i}$ (respectively of $B_{i}$) is a
proper prefix of every word in $B_{i}$ (respectively in $A_{i}$). Note that if
$A_{i}$ and $B_{i}$ are not coinciding singletons, then exactly one of them is
a dominating block.\bigskip

\textbf{Characterization of the Imbalance Order by Comparing Indices\ }%
\textit{For lexicographically enumerated maximal instantaneous codes }
\[
A=a_{1},a_{2},...\text{ \ \ \ \ \ \ \ }B=b_{1},b_{2},...
\]
\textit{of the same size, we have }$A\trianglelefteq B$ \textit{in the
imbalance order (}$A$ \textit{is more balanced than or equal to }$B$\textit{)
if and only if whenever }$a_{i\text{ }}$\textit{and }$b_{j}$ \textit{are
comparable codewords in the prefix order, for their indices we have }$i\geq
j.$

\bigskip

\textbf{Proof} \ Suppose that $A\trianglelefteq B$ and for some
codewords\textit{ }$a_{i\text{ }}$and\textit{ }$b_{j}$ comparable in the
prefix order we have\textit{ }$i<j.$ We shall derive a contradiction. Consider
the interval decompositions for the two codes, as in (2). Due to the
comparability of the two codewords,\ they belong to corresponding intervals,
\textit{i.e.} there is an index $k$ such that $\ a_{i}\in A_{k}$ and $b_{j}\in
B_{k}$. If $A_{k}$ is a singleton, $A_{k}=\left\{  a_{i}\right\}  $, then%
\[
K(a_{1},...,a_{j})>K(a_{1},...,a_{i})=K(A_{1}\cup...\cup A_{k})=K(B_{1}%
\cup...\cup B_{k})\geq K(b_{1},...,b_{j})
\]
contradicts majorization. If $B_{k}$ is a singleton, $B_{k}=\left\{
b_{j}\right\}  $, then majorization is contradicted by
\[
K(a_{1},...,a_{i})>K(A_{1}\cup...\cup A_{k-1})=K(B_{1}\cup...\cup
B_{k-1})=K(b_{1},...,b_{j-1})\geq K(b_{1},...,b_{i})
\]

Suppose conversely that whenever $a_{i\text{ }}$and\textit{ }$b_{j}$ are
comparable codewords in the prefix order, for their indices we have\textit{
}$i\geq j,$ but majorization fails for some index $m$,%
\[
K(a_{1},...,a_{m})>K(b_{1},...,b_{m})\text{
\ \ \ \ \ \ \ \ \ \ \ \ \ \ \ \ \ \ \ \ \ \ \ \ \ \ \ \ \ \ \ \ \ \ \ \ \ \ \ (3)}%
\]
Let the indices $k,l$ be determined by $a_{m}\in A_{k},b_{m}\in B_{l}$. If
$A_{k}$ is a singleton, $A_{k}=\left\{  a_{m}\right\}  $, then (3) requires
that the lex last word $b_{u}$ in $B_{k}$ lexicographically follow $b_{m}.$
Then $m<u$ contradicts the comparability of $a_{m}$ and $b_{u}$. If $B_{l}$ is
a singleton, $B_{l}=\left\{  b_{m}\right\}  $, then (3) requires that $a_{m}$
lexicographically follow all words\ in $A_{l}$. Let $a_{i}$ be any word in
$A_{l}$. Now $i<m$ contradicts the comparability of $a_{i}$ and $b_{m}$.
\ \ \ \ \ \ \ \ \ $\square$

\medskip

\textbf{Reduction Lemma }\textit{If }$B,C$\textit{ are two distinct
lexicographically monotone maximal instantaneous codes of the same size, then
there is a lex monotone maximal instantaneous code }$D$\textit{ that is
(strictly) more balanced than at least one of }$B$\textit{ or }$C$\textit{,
and such that in the imbalance order }$B,C,D$\textit{ have the same common
lower bounds as }$B,C$\textit{.}

\bigskip

\textbf{Proof }Let the given codes be enumerated in lex order as
\[
B=b_{1},b_{2},...\text{ \ \ \ \ \ \ \ }C=c_{1},c_{2},...\text{
\ \ \ \ \ \ \ \ \ \ \ \ \ \ \ \ \ \ \ \ \ \ \ \ \ \ \ \ \ \ \ \ \ \ \ \ \ \ \ \ \ \ \ \ \ \ \ \ \ \ \ \ \ (4)}%
\]
\ \

Since $B$ and $C$ are distinct, there must exist elements $b\in B$ and $c\in
C$ such that

\textit{(i)} $c$ is a proper prefix of some element of $B$

\textit{(ii)} $b$ is a proper prefix of some element of $C.$

Without loss of generality we can assume that the first such $c$
lexicographically precedes the first such $b$. Let $k$ denote the index in $B$
of the lexicographically first element $b$ satisfying condition (ii). With $k$
thus fixed, let $m$ denote the index in $C$\ of the lexicographically last
element $c$ among those elements of $C$ which lexicographically precede
$b_{k}$ and satisfy condition \textit{(i)}. Thus a word $c_{m}$ in $C$ has
also been chosen, and it is easy to see that $m<k$.

With reference to the terminology of decompositions according to the Interval
Decomposition Lemma for Two Codes, $b_{k}$ is the sole element of the first
dominating block of $B$ with respect to $C$, and $c_{m}$ is the sole element
of the last block of $C$ that is dominating and precedes all non-singleton
blocks of $C$.

The code $D$ is now constructed as follows. It is obtained from $C$ by a
single balancing operation, in the sense of $\left[  SPPR\right]  $, chosen to
take into account the relationship of $C$ with $B$, and the choice of the
codewords $c_{m}$ and $b_{k}$. Referring to the indexed enumeration of $C$ in
lex order appearing in (4), let $c_{n},c_{n+1}$ be the first two among the
elements of $C$ admitting $b_{k}$ as a prefix which have the same length. It
is not difficult to see that $m<n$ and $c_{n},c_{n+1}$ are twin sons (one
letter extensions) of some word $w$, $\ c_{n}=w0,c_{n+1}=w1.$

Let $D=(C\setminus\{c_{m},c_{n},c_{n+1}\})$ $\cup$ $\{w,c_{m}0,c_{m}1\}$.
Obviously in the imbalance order $D\trianglelefteq C$.

\bigskip

We claim that if an arbitrary lex monotone maximal instantaneous code $A$ with
lex enumerated codewords $a_{1},a_{2},...$ is more balanced than (or equal to)
$B$ and $C,$ then it is also more balanced than or equal to $D.$ This will
show that the statement of the Lemma holds.

\bigskip

The elements of $D$, enumerated as a sequence of words in lex order as
$d_{1},d_{2},...,$ are partitioned into five consecutive subsequences:

$d_{1}=c_{1\text{ }},$ $...,$ $d_{m-1}=c_{m-1}$
\ \ \ \ \ \ \ \ \ \ \ \ \ \ \ \ \ \ \ \ \ \ \ \ (empty subsequence if $m=1$)

$d_{m}=c_{m}0,$ $d_{m+1}=c_{m}1$

$d_{m+2}=c_{m+1}$,......, $d_{n}=c_{n-1}$ \ \ \ \ \ \ \ \ \ \ \ \ \ \ \ \ \ \ \ \ \ \ \

\ \ \ \ \ \ \ \ \ \ \ \ \ \ \ \ \ \ \ \ \ (empty subsequence if $m+1=n$,
non-empty if $m+1<n$)

$d_{n+1}=w$

$d_{n+2}=c_{n+2\text{ }},d_{n+3}=c_{n+3}$ ,$....$ \ \ \ \ \ \ \ \ \ \ \ \ \ \ \ \ \ \ \ \

\ \ \ \ \ \ \ \ \ \ \ \ \ \ (empty subsequence if $n+1$ is the common size of
the codes)

\bigskip

The subsequence $d_{m+2},...,d_{n}$ in turn consists of two (possibly empty)
consecutive subsequences: the first of these consists of elements also
belonging to $B$, and the second consists of elements that have $b_{k}$ as a
proper prefix. In the first subsequence the index in $B$ of any element is
(strictly) larger than its index in $C$. In the second subsequence the last
symbol of each element is $0$ (i.e. it is not of the form $d_{h}=v1$, for
otherwise $v0$ would have to be also in this second subsequence, contradicting
the definition of $n$).

\bigskip

In view of the Characterization of the Imbalance Order by Comparing Indices,
we need to verify that if some codeword $a$ in $A$ is comparable in the prefix
order to a codeword in $D$, i.e. to some $d_{j}$ having index $j$ in the lex
enumeration of $D,$ then the index of $a$ in $A$ is at least $j.$ This is
obvious for $j<m$ and $n+1<j$, since $A\trianglelefteq B,C$. In the following
examination of the remaining cases \textit{comparability} will always refer to
comparability in the prefix order of the tree of words.

For $j=m,$ if an element $a_{i}$ of $A$ is comparable to $d_{m}=c_{m}0$, then
it is also comparable to its prefix $c_{m}$. The assumption $A\trianglelefteq
C$ then implies $i\geq m$.

For $j=m+1$, if an element $a$ of $A$ is comparable to $d_{m+1}=c_{m}1$, it is
of course also comparable to $c_{m}$, and we claim that it is comparable as
well to some member $b_{h}$ of $B$ with index $h>m.$ Note that both
$d_{m}=c_{m}0$ and \ $d_{m+1}=c_{m}1$\ must be the prefixes of words in $B$,
and no two of the codewords in $C$ can be prefixes of the same word in $B$
(while each one of $c_{1},...,c_{m-1}$ is the prefix of at least one word in
$B$ and $c_{m}$ is a prefix of at least two). From this we can conclude that
$d_{m+1}$ must be the prefix of some $b_{h}$ in $B$, and all such elements of
$B$ have index $h>m.$ Now it follows that the element $a$ of $A$ is comparable
to at least one such $b_{h}$ with index $h>m$. But then, as $A\trianglelefteq
B$ in the imbalance order, the index of $a$ in $A$ is at least $h\geqslant
m+1$.

In the interval $m+2\leq j\leq n,$ if it is not empty, let $j$ be the smallest
index such that for some $i<j$ the elements $a_{i}$ and $d_{j}$ are comparable
- we shall derive a contradiction. For $j$ thus fixed, let $i$ be as small as
possible. If $d_{j}=c_{j-1}$ belongs to $B$, then its index in $B$ is
(strictly) greater than $j-1$, thus by $A\trianglelefteq B$ in the imbalance
order it could not be comparable to $a_{i}$. Therefore $d_{j}=c_{j-1}$ has
$b_{k}$ as a proper prefix. Also its length is (strictly) larger than that of
$c_{j-2}$. This implies that the last symbol of $d_{j}$ must be $0$. But now,
since the last symbol of $d_{j}$ is $0$, if the word $a_{i}$ were a proper
prefix of $d_{j}=c_{j-1}$, then it would be a proper prefix of $c_{j}$ also,
implying $i\geq j$, which is contrary to assumption. Thus $d_{j}$ is a prefix
of $a_{i}$ and $d_{j}$ is the sole element of $D$ comparable to $a_{i}.$ Now
$d_{j-1}$ is comparable to one or more elements $a_{r}$ of $A$, and all such
indices $r$ must be at least $j-1$ by the minimality assumption on $j.$ But as
$a_{i}$ cannot be comparable with $d_{j-1}$, it must come later in the lex
order on $A$ than all the elements $a_{r}$\ of $A$ comparable with $d_{j-1}$,
i.e. $r<i$ for all such $r.$ Therefore $i>r\geq$ $j-1$ and thus $i$ is at
least $j$, a contradiction.

The argument is similar for $j=n+1$. If the element $a_{i}$ of $A$ comparable
with $d_{n+1}=w$ is comparable with $c_{n+1}=w1$ then we are done. Else
$a_{i}$ must have $c_{n}=w0$ and $w=d_{n+1}$ as a prefix and cannot be
comparable with $d_{n}$. Therefore, $i>r\geq j-1$ and thus $a_{i}$ must come
later in the lex order on $A$ than all the elements of $A$ comparable with
$d_{n}$, i.e. $r<i$ for all such $r$. But we already know that the indices in
$A$ of these latter elements $a_{r}$ are at least $n,$ forcing $i\geqslant
n+1.$

\bigskip

We have thus shown that in the imbalance order all common lower bounds $A$ of
$B$ and $C$\ are also lower bounds of the code $D$ constructed from these
latter two, completing the proof of the Lemma and thus providing an
alternative proof of the Lattice Property. \ \ \ \ \ \ \ $\square$

\bigskip

\textit{Remark. }Repeated application of the construction of $D$ in the proof
of the Reduction Lemma provides an algorithm for constructing the meet of any
two codes $B$ and $C$ in the imbalance lattice. (The repetition is to be
applied\ to the reduced pair of codes gotten by replacing $B$ or $C$ by $D$,
according to whether $D$ is more balanced than $B$ or $C.$) \ As simple
examples with incomparable $B$ and $C$, we can take, using the lattice
diagrams on p. 7. of \ $\left[  SPPR\right]  ,$\ with path-length sequence
representation of codes of size $7,$ $\bigskip$

$B$ $=$ $2$ $2$ $2$ $3$ $4$ $5$ $5$ \ \ \ \ $C$ $=$ $1$ $3$ $3$ $4$ $4$ $4$
$4$\ \ \ \ $D$ $=$ $2$ $2$ $2$ $4$ $4$ $4$ $4$ \ \ \ \

(where $D$ is in fact the meet of $B$ and $C$)\bigskip

\noindent or we can take as example of codes of size $9$

$\bigskip$

$B$ $=$ $233\ 333$ $455$ \ \ \ \ \ \ \ \ $C$ $=$ $144$ $444\ 444$
\ \ \ \ \ \ $D$ $=$ $223$ $444\ 444$\ \ \ \ \

(where $D$ is still less balanced than the meet of $B$ and $C$).

\bigskip

\bigskip

\textbf{Acknowledgements.}

Part of this work has been co-funded by Marie Curie Actions and supported by
the National Development Agency (NDA) of Hungary and the Hungarian Scientific
Research Fund (OTKA, contract number 84593), within a project hosted by the
University of Miskolc, Department of Analysis.

\vskip 1cm

\includegraphics[height=0.9cm]{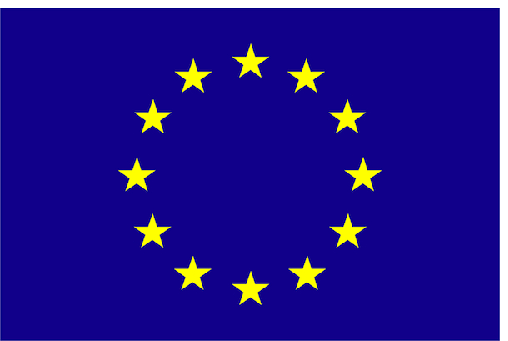} \quad
\includegraphics[height=1.2cm]{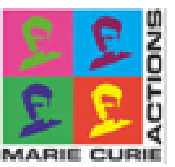} \quad
\includegraphics[height=1.2cm]{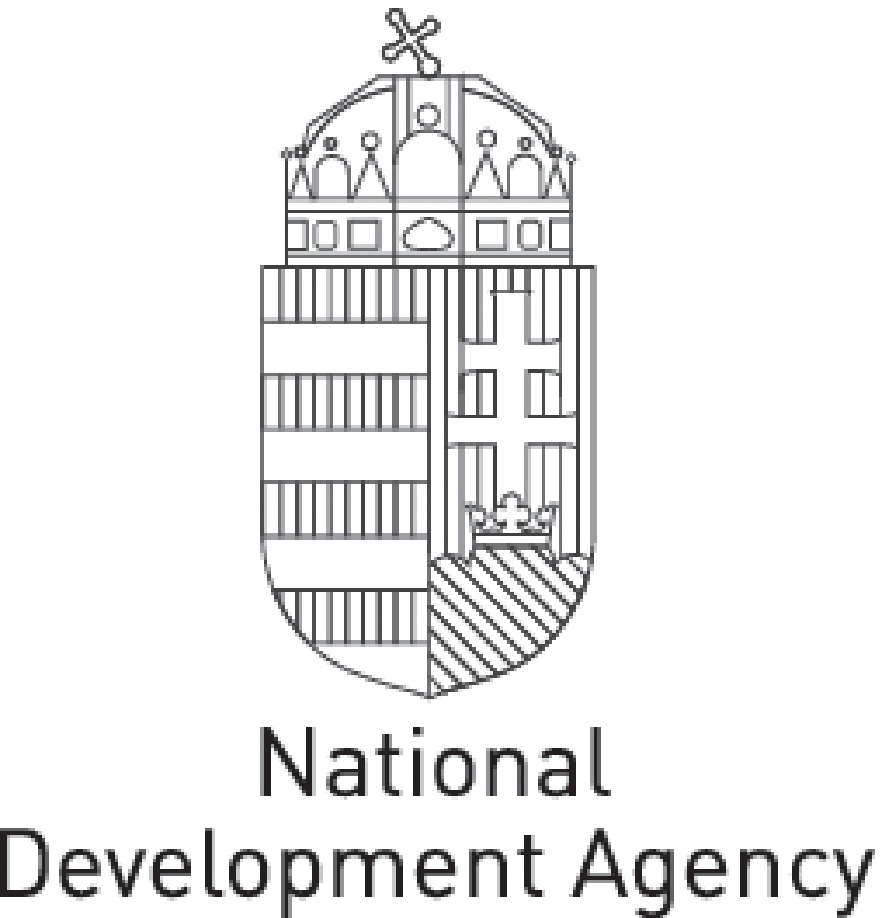} \quad
\includegraphics[height=1cm]{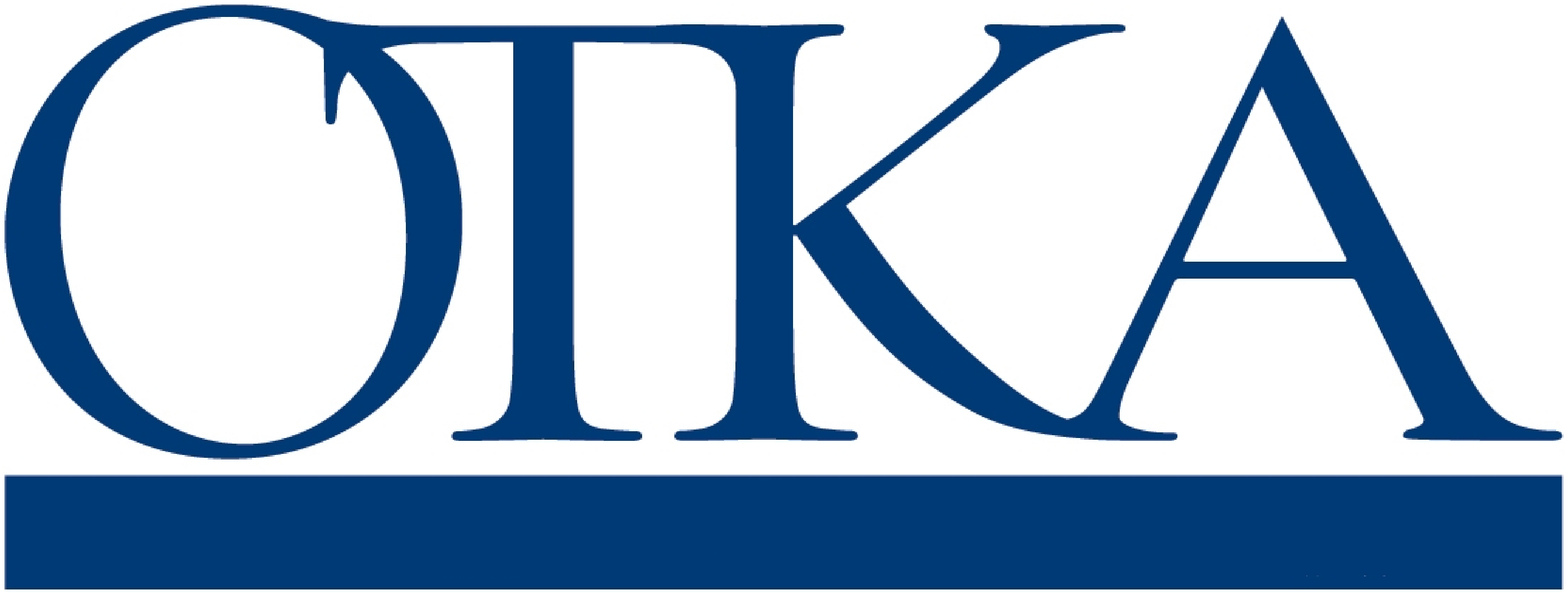}

\vskip 1cm

\textbf{References}

\bigskip

$\left[  EHP\right]  $ C. Elsholtz, C. Heuberger, H. Prodinger, The number of
Huffman codes, compact trees, and sums of unit fractions, \textit{IEEE Trans.
Information Theory, }59 (2) 2013, 1065-1075

$\bigskip$

$\left[  FR\right]  $ \ S.\ Foldes, S. Radeleczki, On the imbalance lattice of
path-length sequences of binary trees, ArXiv 2013 (https://arxiv.org/abs/1307.0161)

\bigskip

$\left[  FS\right]  $ \ S. Foldes, N.M. Singhi, On instantaneous codes,
Journal of Combinatorics, Information \& System Sciences 31 (2006) 317--326

\bigskip

$\left[  K\right]  $ \ L.G. Kraft, \textit{A Device for Quantizing, Grouping,
and Coding Amplitude Modulated Pulses}, Q.S. Thesis, MIT 1949

$\bigskip$

$\left[  OKPS\right]  $ \ M. O'Keefe, H. Pajoohesh, M. Schellekens, Decision
trees of algorithms and a semivaluation to measure their distance,
\textit{Electr. Notes Comput. Sc. }161 (2006) 175-183

$\bigskip$

$\left[  P\right]  $ H. Pajoohesh, Topological and categirical properties of
of binary trees, \textit{Applied Gen. Topology} 9 (1) (2008) 1-14

\bigskip

$\left[  SPPR\right]  $ \ D. Stott Parker, Prasad Ram, The Construction of
Huffman Codes is a Submodular ("Convex") Optimization Problem Over a Lattice
of Binary Trees. \textit{SIAM J. Comput.} 28(5) 1875-1905 (1999)

\bigskip

\bigskip

\textbf{Authors' addresses:}

\bigskip

S. Foldes

\texttt{stephan@renyi.mta.hu}

\medskip

D.S. Parker

UCLA Computer Science Department

\texttt{stott@cs.ucla.edu}

\medskip

S. Radeleczki

University of Miskolc, Institute of Mathematics

\texttt{matradi@uni-miskolc.hu}
\end{document}